\documentclass[12pt,reqno]{amsart}

\usepackage[colorlinks=true,linkcolor=blue,citecolor=blue,urlcolor=blue]{hyperref}
\usepackage{array}
\usepackage{amssymb,latexsym}
\usepackage{bm,bbding}
\usepackage{cleveref}
\usepackage{xcolor}
\usepackage{pstricks-add}
\usepackage[all]{xy}

\newtheorem{theorem}{Theorem}[section]

\newtheorem{corollary}[theorem]{Corollary}

\numberwithin{equation}{section}
\theoremstyle{definition}
\newtheorem{remark}[theorem]{Remark}

\newcommand{\rr}{\mathbb{R}}

\newcommand{\zz}{\mathbb{Z}}

\newcommand{\td}{{\rm td}}
\newcommand{\Td}{{\rm Td}}

\begin{document}

\title{Todd genus and $A_k$-genus of unitary $S^1$-manifolds}\thanks{This project was supported by the Natural Science Foundation of Tianjin City of China, No. 19JCYBJC30300.}

\author{Jianbo Wang, Zhiwang Yu, Yuyu Wang}
\address{(Jianbo Wang) School of Mathematics, Tianjin University, Tianjin 300350, China}
\email{wjianbo@tju.edu.cn}
\address{(Zhiwang Yu) School of Mathematics, Tianjin University, Tianjin 300350, China}
\email{yzhwang@tju.edu.cn}
\address{(Yuyu Wang) College of Mathematical Science, Tianjin Normal University, Tianjin 300387, China}
\email{wdoubleyu@aliyun.com}

\begin{abstract} 
Assume that $M$ is a compact connected unitary $2n$-dimensional manifold and admits a non-trivial circle action preserving the given complex structure. If the first Chern class of $M$ equals to $k_0x$ for a certain $2{\rm nd}$ integral cohomology class $x$ with $|k_0|\geqslant n+2$, and its first integral cohomology group is zero,  this short paper shows that  the Todd genus and $A_k$-genus of $M$ vanish.
\end{abstract}

\date{\today}
\maketitle

\section{Introduction}
A smooth manifold $M$ is called a \textit{unitary manifold} (sometimes being called {\it stably almost complex manifold} or {\it weakly almost complex manifold}), if the tangent bundle of $M$ admits a stably complex structure. Namely, there exists a bundle map
\[
J: TM \oplus \underline{\rr}^{l} \longrightarrow TM \oplus  \underline{\rr}^{l} \\
\]
such that $J^{2}=-1$, where $ \underline{\rr}^l$ denotes the trivial real $l$-plane bundle over $M$ for some $l$. A stable complex structure induces an orientation, obtained as the “difference” of the complex orientation on $TM\oplus\underline{\rr}^l$ and the standard orientation on $\underline{\rr}^l$. Hereafter a unitary manifold is always oriented in such a way. If $S^1$ acts smoothly on a unitary (respectively almost complex) manifold $M$ and if the differential of each element of $S^1$ preserves the given complex vector bundle structure then $M$ will be called \textit{unitary} (respectively \textit{almost complex}) \textit{$S^1$-manifold}.


Given some assumptions, there are many results about the existence of $S^1$ actions on manifolds with the vanishing of Todd genus, $A_k$-genus or $\hat{A}$-genus. Let's list some related results.
\begin{itemize}
\item Hattori proved that (\cite[Proposition 3.21]{Hattori1984}): for a unitary $S^{1}$-manifold $M$ having only isolated fixed points, if the first Chern class is a torsion element then the Todd genus $\Td(M)$ is zero. 
\item Kri$\check{c}$ever showed that (\cite[Theorem 2.2]{Kricever1976}), for a unitary $S^1$-manifold $M$,  if the first Chern class $c_1(M)$ is divisible by $k$, the $A_k$-genera of $M$ vanish, i.e., $A_k(M)=0, k\geqslant 2$.
\item
Atiyah and Hirzebruch \cite{AtHi70} gave the result that if $M$ is a connected $2n$-dimensional spin manifold and $S^1$ acts non-trivially on $M$, then the $\hat{A}$-genus $\hat{A}(M)$ vanishes.
\item Using Hattori's result (\cite[Proposition 3.21]{Hattori1984}) and  the celebrated theorem of Atiyah and Hirzebruch in \cite{AtHi70},   Li (\cite[Proposition 4.1]{LiPing2011}) showed that: for an almost complex manifold admitting a non-trivial $S^{1}$-action, if the first Chern class is a torsion element and is zero under $\mod 2$ reduction, then the Todd genus is zero.
\item
In \cite[Theorem 1.1]{FangRong2000}, Fang and Rong showed that:
For a compact $2n$-manifold $M$ of finite fundamental group, if $M$ admits a fixed point free $S^{1}$-action, then for all $c\in H^{2}(M;\zz)$ and $0 \leqslant m \leqslant\left[\frac{n}{2}\right]$, $\left(c^{n-2 m} \cdot L_{m}\right)[M]=0$, where $L= L_0+L_1+\cdots+L_{\left[\frac{n}{2}\right]}$ is the Hirzebruch $L$-polynomial of $M$.
In particular, $\left(c^{n}\right)[M]=0$. 
\item A generalization of Fang-Rong's result is the \cite[Example 3.2]{Li2012} proved by  Li: Let $M^{2 m}$ be a manifold with $b_{1}(M)=0 .$ If $M$ admits a fixed point free $S^{1}$-action, then for all $c \in H^{2}(M ; \mathbb{Z})$ and $0 \leqslant j \leqslant\left[\frac{m}{2}\right]$,  $\left(c^{m-2 j} \cdot L_{j}\right)[M]=0$,
where $L_{j}$ is any polynomial of degree $j$ in the subalgebra of $H^{*}(M)$ generated by Pontrjagin classes $p_{i}(M)\left(\operatorname{deg}\left(p_{i}\right)=i\right)$. In particular, $\left(c^{n}\right)[M]=0$.
\end{itemize}
As a supplement to the above mentioned results,
we  prove the following theorem, which is the main theorem of this paper.

\begin{theorem}\label{Toddgenusvanish}
Let $M$ be a compact connected unitary manifold of dimension $2n>2$ with $H^1(M;\zz)=0$. Suppose that $M$ admits a non-trivial $S^1$-action which preserves the complex structure of the stable tangent bundle of $M$. If the first Chern class $c_1(M)=k_0x$ for some $x\in H^2(M;\zz)$ and $|k_0|\geqslant n+2$, then

(1) $x^{n-2j}\hat{A}_{j}=0$, where $\hat{A}_j:=\hat{A}_j(p_1,\dots,p_j)$ are the $\hat{A}$ polynomials, $0\leqslant j\leqslant \left[\frac{n}{2}\right]$.

(2) $x^{n-i}T_{i}=0$, where $T_i:=T_i(c_1,\dots,c_i)$ are the Todd polynomials, $0\leqslant i\leqslant n$. In particular, the Todd genus $\Td(M)=T_n[M]=0$.

(3) The $A_k$-genus $A_k(M)=0$, $k\geqslant 2$.
\end{theorem}
The proof of \Cref{Toddgenusvanish} is mainly based on the following theorem of Hattori \cite[Thoerem 3.13]{Hattori1978}.
\begin{theorem}[Hattori]\label{Hattori-thm3.13}
Let $M$ be a closed connected weakly almost complex manifold of dimension $2n>2$ with $H^1(X;\zz)=0$. Suppose that $M$ admits a non-trivial $S^1$-action which preserves the given complex structure of the stable tangent bundle of $M$. If $c_1(M)=k_0x, x\in H^2(M;\zz)$, then we have
\[
\left\{\exp\left(\frac{kx}{2}\right)\cdot \hat{A}(TM)\right\}[M]=0
\]
for each integer $k$ such that $k\equiv k_0 \pmod 2$ and $|k|<|k_0|$.
\end{theorem}
\begin{remark}
(1) In \cite[Corollary 4.2]{Sabatini2017}, for an almost complex $S^1$-manifold $M$ of dimension $2n$ with only isolated fixed points, Sabatini showed that if $c_1(M)=k_0x$ for some $x\in H^2(M;\zz)$ and $k_0\geqslant n+2$, then $c_1^{n-k}T_k=0$, $0\leqslant k\leqslant n$. Although Sabatini's conclusion is similar to our \Cref{Toddgenusvanish} (2), the condition of isolated fixed points is not necessary in \Cref{Toddgenusvanish}.

(2) The vanishing of $A_k$-genera in \Cref{Toddgenusvanish} (3) hold for any $k\geqslant 2$, and the vanishing of $A_k$-genera in \cite[Theorem 2.2]{Kricever1976} hold for those $k$ satisfying that the first Chern class is divisible by $k$.
\end{remark}

\section{Preliminaries}

Firstly, let's briefly introduce the Todd genus, $A_k$-genus and $\hat{A}$-genus.

Let $M$ be a connected, closed $2n$-dimensional unitary manifold. Since the Chern classes are stable invariants, which means that they are unchanged if we add a trivial (complex) bundle, the structure one needs for defining the Todd genus is a stable almost complex structure. For abbreviation, the Chern classes of the unitary manifold $M$ are denoted as 
\(
c_i:=c_i(M):=c_i(TM\oplus\underline{\rr}^l). 
\)
The Todd genus $\Td(M)$ of $M$ is the evaluation of the Todd polynomial  $T_n(c_1,c_2,\dots,c_n)$ in the Chern classes on the fundamental class $[M]$ of $M$. The multiplicative sequence $\{T_k(c_1,\dots,c_k)\}$ is associated to the power series $\displaystyle
Q(x)=\frac{x}{1-\exp(-x)}$ (\cite[\S 1.7]{Hi1978}).
Define the Todd class of $M$ by
\[
\td(TM\oplus\underline{\rr}^l)=\sum_{k=0}^\infty{T_k(c_1,\dots,c_k)},
\]
then the Todd genus of $M$ is  the evaluation
 \[
 \Td(M)=\td(TM\oplus\underline{\rr}^l)[M]=T_n(c_1,\dots,c_n)[M].
 \]
Associated with the characteristic power series $Q(x)=\dfrac{kx\exp(x)}{\exp(kx)-1}$, $k\geqslant 2$, there is a multiplicative $A_k$-class  of $M$  (\cite{Kricever1976}), which we denote it as $A_k(TM\oplus\underline{\rr}^l)$, and the $A_k$-genus $A_k(M)$ is the evaluation
 \[
 A_k(M)=A_k(TM\oplus\underline{\rr}^l)[M].
 \]
 Note that, if $k=1$ is allowed in \cite{Kricever1976}, $A_1(M)$ is exactly the Todd genus $\Td(M)$.
If we consider a formal factorization of the total Chern class as $c(TM\oplus\underline{\rr}^l)=(1+x_1)\cdots(1+x_n)$, then by \cite[\S1.8]{HiBeJu1992},
\begin{align*}
\Td(M) & =\left(\prod_{i=1}^n \dfrac{x_i}{1-\exp(-x_i)}\right)[M],\\
A_k(M) & =\left(\prod_{i=1}^n\dfrac{kx_i\exp(x_i)}{\exp(kx_i)-1}\right)[M].
\end{align*}

 For a compact oriented differentiable $2n$-manifold  $M$, Hirzebruch (\cite[page 197]{Hi1978}) defined the $\hat{A}$-sequence of $M$ as a certain polynomials in the Pontrjagin classes of $M$. More concretely, the even power series
\[
Q(x)=\frac{\frac{1}{2}x}{\sinh \left(\frac{1}{2}x\right)}=1-\frac{x^2}{24}+\frac{7x^4}{5760}-\frac{31x^6}{967680}+\frac{127x^8}{154828800}+\cdots
\]
defines a multiplicative sequence $\left\{\hat{A}_{j}(p_{1}, \ldots, p_{j})\right\}$ (for abbreviation, $p_j:=p_j(TM)$), where $\hat{A}_j(p_{1}, \ldots, p_{j})$ is a rational homogeneous polynomial of degree $4j$  in the Pontrjagin classes.
The $\hat{A}$-class $\hat{A}(TM)$ is defined as follows:
\begin{equation*}
\hat{A}(TM)=\sum_{j=0}^\infty{\hat{A}_j(p_1,\dots,p_j)}.
\end{equation*}
The $\hat{A}$-genus of $M$ is the evaluation
\[
\hat{A}(M)=\hat{A}(TM)[M].
\]
Let $p(TM)=(1+x_1^2)\cdots(1+x_n^2)$ be the formal factorization of the total Pontrjagin class, then by \cite[\S1.6]{HiBeJu1992},
\[
\hat{A}(M)=\left(\prod_{i=1}^n\frac{\frac{1}{2}x_i}{\sinh \left(\frac{1}{2}x_i\right)}\right)[M].
\]

\section{Proof of  the main \Cref{Toddgenusvanish}}\label{sec-vanish}

\begin{proof}[Proof of \Cref{Toddgenusvanish}]
(1) Under the assumptions of \Cref{Toddgenusvanish}, by \Cref{Hattori-thm3.13},
\begin{equation}\label{HattoriTHM}
\left\{\exp\left(\frac{kx}{2}\right)\cdot \hat{A}(TM)\right\}[M]=0
\end{equation}
for each integer $k$ such that  $|k|<|k_0|, k \equiv k_0 \mod 2$. Let $u^{(2n)}$ be the term of degree $2n$ of an element $u\in\bigoplus\limits_{k=0}^{2n} H^k(M;\zz)$. By the Poincar\'e duality, from \eqref{HattoriTHM} we can obtain

\begin{equation*}
\left(\exp\left(\frac{kx}{2}\right)\cdot \hat{A}(TM)\right)^{(2n)}=0\in H^{2n}(M;\zz).
\end{equation*}


{\bf Case 1:} $n$ is even.

We have
\begin{align*}
\left(\exp\left(\frac{kx}{2}\right)\cdot\hat{A}(TM)\right)^{(2n)}  = & \frac{1}{n!}\left(\frac{kx}{2}\right)^n+\frac{1}{(n-2)!}\left(\frac{kx}{2}\right)^{n-2}\hat{A}_1+\cdots\\
& +\frac{1}{2!}\left(\frac{kx}{2}\right)^{2}\hat{A}_{\frac{n}{2}-1}+\hat{A}_{\frac{n}{2}}. 
\end{align*}
When $|k_0|\geqslant n+2$, there exist integers $k_1,k_2,\dots,k_{\frac{n}{2}+1}$ such that
\begin{align}
& |k_1|<|k_2|<\cdots <|k_{\frac{n}{2}}|<|k_{\frac{n}{2}+1}|<|k_0|,\nonumber\\
& k_1\equiv k_2\equiv\cdots\equiv k_{\frac{n}{2}}\equiv k_{\frac{n}{2}+1}\equiv k_0 \mod 2, \nonumber\\
\text{and~~}& \left(\exp\left(\frac{k_ix}{2}\right)\cdot \hat{A}(TM)\right)^{(2n)}=0, ~i=1,2,\dots,\frac{n}{2}+1. \label{expAhateven}
\end{align}
By \eqref{expAhateven}, we get the following equation of matrices:
\begin{equation*}
\begin{bmatrix}
k_1^n&k_1^{n-2}&\cdots&k_1^2&1\\
k_2^n&k_2^{n-2}&\cdots&k_2^2&1\\
\vdots&\vdots&\cdots&\vdots&\vdots\\
k_{\frac{n}{2}}^n&k_{\frac{n}{2}}^{n-2}&\cdots&k_{\frac{n}{2}}^2&1\\
k_{\frac{n}{2}+1}^n&k_{\frac{n}{2}+1}^{n-2}&\cdots&k_{\frac{n}{2}+1}^2&1
\end{bmatrix}
\begin{bmatrix}
\frac{1}{n!\cdot 2^n}x^n\\
\frac{1}{(n-2)!\cdot 2^{n-2}}x^{n-2}\hat{A}_1\\
\vdots\\
\frac{1}{2!\cdot 2^2}x^{2}\hat{A}_{\frac{n}{2}-1}\\
\hat{A}_{\frac{n}{2}}
\end{bmatrix}=0.
\end{equation*}
Since $k_1,k_2,\dots,k_{\frac{n}{2}+1}$ are mutually distinct, which means the determinant of left matrix in the above equation is that of the Vandermonde and this matrix is invertible, so we get
\begin{equation}\label{evenxAhat}
x^n=0, ~x^{n-2}\hat{A}_1=0, \dots, ~x^{2}\hat{A}_{\frac{n}{2}-1}=0,  ~\hat{A}_{\frac{n}{2}}=0.
\end{equation}

{\bf Case 2:} $n$ is odd.

We have
\begin{align*}
\left(\exp\left(\frac{kx}{2}\right)\cdot\hat{A}(TM)\right)^{(2n)}  = & \frac{1}{n!}\left(\frac{kx}{2}\right)^n+\frac{1}{(n-2)!}\left(\frac{kx}{2}\right)^{n-2}\hat{A}_1+\cdots\\
& +\frac{1}{3!}\left(\frac{kx}{2}\right)^{3}\hat{A}_{\frac{n-3}{2}}+\frac{1}{1!}\left(\frac{kx}{2}\right)\hat{A}_{\frac{n-1}{2}}. 
\end{align*}
When $|k_0|\geqslant n+2$, there exist integers $k_1,k_2,\dots,k_{\frac{n+1}{2}}$ such that
\begin{align}
& 0<|k_1|<|k_2|<\cdots <|k_{\frac{n-1}{2}}|<|k_{\frac{n+1}{2}}|<|k_0|, \nonumber\\
& k_1\equiv k_2\equiv\cdots\equiv k_{\frac{n-1}{2}}\equiv k_{\frac{n+1}{2}}\equiv k_0 \mod 2, \nonumber\\
\text{and~~}& \left(\exp\left(\frac{k_ix}{2}\right)\cdot \hat{A}(TM)\right)^{(2n)}=0, ~i=1,2,\dots,\frac{n+1}{2}. \label{expAhatodd}
\end{align}
By \eqref{expAhatodd}, we get the following equation of matrices:
\begin{equation*}
\begin{bmatrix}
k_1^n&k_1^{n-2}&\cdots&k_1^3&k_1\\
k_2^n&k_2^{n-2}&\cdots&k_2^3&k_2\\
\vdots&\vdots&\cdots&\vdots&\vdots\\
k_{\frac{n-1}{2}}^n&k_{\frac{n-1}{2}}^{n-2}&\cdots&k_{\frac{n-1}{2}}^3&k_{\frac{n-1}{2}}\\
k_{\frac{n+1}{2}}^n&k_{\frac{n+1}{2}}^{n-2}&\cdots&k_{\frac{n+1}{2}}^3&k_{\frac{n+1}{2}}
\end{bmatrix}
\begin{bmatrix}
\frac{1}{n!\cdot 2^n}x^n\\
\frac{1}{(n-2)!\cdot 2^{n-2}}x^{n-2}\hat{A}_1\\
\vdots\\
\frac{1}{3!\cdot 2^3}x^{3}\hat{A}_{\frac{n-3}{2}}\\
\frac{1}{2}x\hat{A}_{\frac{n-1}{2}}
\end{bmatrix}=0.
\end{equation*}
Since the determinant of left matrix in the above equation is a multiple of that of the Vandermonde,  we get
\begin{equation}\label{oddxAhat}
x^n=0, ~x^{n-2}\hat{A}_1=0, \dots, ~x^{3}\hat{A}_{\frac{n-3}{2}}=0,  ~x\hat{A}_{\frac{n-1}{2}}=0.
\end{equation}

(2)
By \cite[(12)~in page 13]{Hi1978} and $A_s=2^{4s}\hat{A}_s$ in \cite[page 197]{Hi1978}, we have
\begin{align*}
T_{k} & =\sum_{r+2s=k} \frac{1}{r!\cdot 2^r} c_{1}^{r} \hat{A}_{s}\\
& =\sum_{s=0}^{[\frac{k}{2}]} \frac{k_0^{k-2s}}{(k-2s)!\cdot 2^{k-2s}} x^{k-2s} \hat{A}_{s},
\end{align*}
where $r$ and $s$ are non-negative integers.
Then for any $0\leqslant k\leqslant n$,
\begin{align*}
x^{n-k}T_k & =x^{n-k}\cdot \sum_{s=0}^{[\frac{k}{2}]} \frac{k_0^{k-2s}}{(k-2s)!\cdot 2^{k-2s}} x^{k-2s} \hat{A}_{s}\\
&=\sum_{s=0}^{[\frac{k}{2}]} \frac{k_0^{k-2s}}{(k-2s)!\cdot 2^{k-2s}} x^{n-2s} \hat{A}_{s}\\
&=0.
\end{align*}

(3) Let 
$$A_{k}(TM\oplus\underline{\rr}^l), ~~A_{\frac{1}{k}}(TM\oplus\underline{\rr}^l)$$ 
be the multiplicative characteristic classes of $M$ associated to the characteristic power series 
\[
\dfrac{kx\exp(x)}{\exp(kx)-1}, ~~\dfrac{x\exp\left(\frac{x}{k}\right)}{\exp(x)-1},~k\geqslant 2,
\]  
and $A_k(M)$, $A_{\frac{1}{k}}(M)$ denote the corresponding genus of $M$ respectively,  then
\begin{equation}\label{AkAk-1}
A_k(M)=k^n\cdot A_{\frac{1}{k}}(M).
\end{equation}
The identity
\[
\frac{x\exp\left(\frac{x}{k}\right)}{\exp(x)-1}=\exp\left(\left(\frac{1}{k}-\frac{1}{2}\right)x\right)\cdot\frac{\frac{x}{2}}{\sinh\left(\frac{x}{2}\right)}
\]
implies that
\[
A_{\frac{1}{k}}(TM\oplus\underline{\rr}^l)=\exp\left(\left(\frac{1}{k}-\frac{1}{2}\right)c_1(M)\right)\cdot\hat{A}(TM).
\]
Thus, as the proof of of \Cref{Toddgenusvanish} (1), using \eqref{evenxAhat} and \eqref{oddxAhat}, we also get $A_{\frac{1}{k}}(TM\oplus\underline{\rr}^l)$ are zero, then the $A_{\frac{1}{k}}$-genera and $A_{k}$-genera of $M$ vanish, $k\geqslant 2$.
\end{proof}
\begin{remark}
In fact,  for $A_2$-genus, the corresponding power series 
\(
\dfrac{2x\cdot \exp{(x)}}{\exp{(2x)}-1}
\)
 is an even power series, then $A_2$-genus is expressible in Pontrjagin numbers and hence defined for an oriented smooth manifold (\cite[\S1.6]{HiBeJu1992}). By \cite[Appendix III]{HiBeJu1992}, 
\[
A_{\frac{1}{2}}(M)=\chi(M,K^{1/2})=\hat{A}(M), 
\] 
where, $\chi(M,K^{1/2})$ is the genus with respect to the characteristic power series $\dfrac{x\exp\left(\frac{x}{2}\right)}{\exp(x)-1}$.
So by \eqref{AkAk-1}, we have 
\[
 A_2(M)=2^n\cdot\hat{A}(M).
\]
\end{remark}

\begin{corollary}
Under the assumptions in \Cref{Toddgenusvanish},  if $c:=c_1(M)\in H^2(M;\zz)$ is a torsion element, then
\begin{align*}
& (c^{n-i}T_i)[M]=0, 0\leqslant i\leqslant n;\\
& (c^{n-2j}\hat{A}_j)[M]=0, 0\leqslant j\leqslant \left[\frac{n}{2}\right];\\
& A_k(M)=0, k\geqslant 2.
\end{align*}
\end{corollary}

\renewcommand\refname{References}

\end{document}